\documentclass[12pt,twoside]{article}

\begin{document}

\centerline{}

\centerline {\Large{\bf Relations on intuitionistic fuzzy soft sets}}

\centerline{}

\newcommand{\mvec}[1]{\mbox{\bfseries\itshape #1}}

\centerline{\bf {Bivas Dinda$^1$ and  T.K. Samanta$^2$}}
\centerline{}

\centerline{$^1$Department of Mathematics,}
\centerline{Mahishamuri Ramkrishna Vidyapith}\centerline{Howrah-711401, West Bengal, India }
\centerline{e-mail: bvsdinda@gmail.com}
\centerline{}
\centerline{$^2$Department of Mathematics}\centerline{ Uluberia
College}\centerline{Howrah, West Bengal, India.}
\centerline{e-mail: mumpu$_{-}$tapas5@yahoo.co.in}
\centerline{}

\newtheorem{Theorem}{\quad Theorem}[section]

\newtheorem{definition}[Theorem]{\quad Definition}

\newtheorem{theorem}[Theorem]{\quad Theorem}

\newtheorem{remark}[Theorem]{\quad Remark}

\newtheorem{corollary}[Theorem]{\quad Corollary}

\newtheorem{note}[Theorem]{\quad Note}

\newtheorem{lemma}[Theorem]{\quad Lemma}

\newtheorem{example}[Theorem]{\quad Example}

\newtheorem{proposition}[Theorem]{\quad Proposition}

\begin{abstract}
\textbf{\emph{In this paper, we present the concept of relations in intuitionistic fuzzy soft set and study some of their properties and also discuss symmetric, transitive and reflexive intuitionistic fuzzy soft relations. }}
\end{abstract}

{\bf Keywords:}  \emph{Intuitionistic fuzzy soft set, union, intersection, cartesian product, intuitionistic fuzzy soft relations.}\\

{\bf 2000 MSC No:} 08A02.

\section{Introduction}
Most of the problems of real life have various uncertainties. Traditional mathematical tools are unable to solve uncertain problems. There are theories viz. theory of probability, theory of evidence, fuzzy set, intuitionistic fuzzy set, vague set for dealing with uncertainties. Thease theories have their own difficulties. The reason for these difficulties is inadequacy of parametrization tool of the theories. Molodtsov proposed the novel concept of soft set theory in his pioneering paper \cite{Molodtsov}. Later on  authors like
Maji et al. \cite{Maji2,Maji,Maji1} have further studied the theory of soft sets and introduced the concepts of fuzzy soft set and intuitionistic fuzzy soft set.\\
In this paper we have introduced the concept of intuitionistic fuzzy soft relation and studied a few of its algebraic properties. We have also discussed symmetric, transitive and reflexive intuitionistic fuzzy soft relations.\\
The organization of this paper is as follows: In section 2, some basic definition and preliminary results are given which will be used in the rest of the paper. In section 3, union and intersection of intuitionistic fuzzy soft sets redefined. In section 4, relation on intuitionistic fuzzy soft set is defined and some of its algebraic properties are studied. In section 5, symmetric, transitive and reflexive intuitionistic fuzzy soft relations are defined and a few their properties are studied.

\section{Preliminaries}
\begin{definition}\cite{Molodtsov}
Let $U$ be an initial universe set and $E$ be the set of parameters. Let $P(U)$ denotes the power set of $U$. A pair $(F,E)$ is called a soft set over $U$ where $F$ is a mapping given by $F:E\rightarrow P(U)$.
\end{definition}

\begin{definition}\cite{Maji2}
Let $U$ be an initial universe set and $E$ be the set of parameters. Let $A \subset E$. A pair $(F,A)$ is called fuzzy soft set over $U$ where $F$
is a mapping given by $F : A\rightarrow I^U  $, where $I^U$ denotes the collection of all fuzzy subsets of $U$.
\end{definition}

\begin{definition}\cite{Maji}
Let $U$ be an initial universe set and $E$ be the set of parameters. Let $IF^U$ denotes the collection of all intuitionistic fuzzy subsets of $U$. Let $A\subseteq E$. A pair $(F,\,A)$ is called intuitionistic fuzzy soft set over $U$, where $F$ is a mapping given by $F:\,A\rightarrow \,IF^U.$
\end{definition}

\begin{example}
Consider the following example:\\
Let $(F,A)$ describes the the character of the students with respect to the given parameters, for finding the best student of an academic year. Let the set of students under consideration is $U=\{s_1,s_2,s_3,s_4\}.$ Let $A\subseteq E$ and $A=\{r="result",\,c="conduct",\,g="games\; and\; sports\; performances"\}$. Let\\
$F(r)=\left\{(s_1,0.8,0.1),(s_2,0.9,0.05),(s_3,0.85,0.1),(s_4,0.75,0.2)\right\}\\$
$F(c)=\left\{(s_1,0.6,0.3),(s_2,0.65,0.2),(s_3,0.7,0.2),(s_4,0.65,0.2)\right\}\\$
$F(g)=\left\{(s_1,0.75,0.2),(s_2,0.5,0.3),(s_3,0.5,0.4),(s_4,0.7,0.2)\right\}$.\\
Then the family $\{F(r),F(c),F(g)\}$ of $IF^U$ is an intuitionistic fuzzy soft set.
\end{example}

\begin{definition}\cite{Maji}
Intrersection of two intuitionistic fuzzy soft sets $(F,A)$ and $(G,B)$ over a common universe $U$ is the intuitionistic fuzzy soft set $(H,C)$ where $C=A\cap B,$ and $\forall \epsilon\in C,\;H(e)=F(e)\cap G(e).$ We write $(F,A)\tilde{\cap}(G,B)=(H,C).$
\end{definition}

\begin{definition}\cite{Maji}
Union of two intuitionistic fuzzy soft sets $(F,A)$ and $(G,B)$ over a common universe $U$ is the intuitionistic fuzzy soft set $(H,C)$ where $C=A\cup B,$ and $\forall \epsilon\in C,$
\[H(e)=F(e),\hspace{3 cm}if\; e\in A-B\]
\[=G(e),\hspace{3 cm}if\;e\in B-A\hspace{-1 cm}\]
\[=F(e)\cup G(e),\hspace{1.75 cm}if\; e\in A\cap B\hspace{-1 cm}\]
We write $(F,A)\tilde{\cup}(G,B)=(H,C).$
\end{definition}

\begin{definition}\cite{Maji}
For two intuitionistic fuzzy soft sets $(F,A)$ and $(G,B)$ over a common universe $U$, we say that $(F,A)$ is an intuitionistic fuzzy soft subset of $(G,B)$ if\\
(i) $A\subset B,$ and\\
(ii) $\forall \epsilon\in A,\;F(\epsilon)$ is an intuitionistic fuzzy subset of $G(\epsilon).\\$
We write $(F,A)\tilde{\subset}(G,B).$
\end{definition}

\begin{definition}
\cite{Schweizer} A binary operation \, $\ast \; : \; [\,0 \; , \;
1\,] \; \times \; [\,0 \; , \; 1\,] \;\, \rightarrow \;\, [\,0
\; , \; 1\,]$ \, is continuous \, $t$ - norm if \,$\ast$\, satisfies
the
following conditions \, $:$ \\
(i)  $\ast$ \, is commutative and associative ,\\
(ii) $\ast$ \, is continuous , \\
(iii) $a \;\ast\;1 \;\,=\;\, a \hspace{1.2cm}
\forall \;\; a \;\; \in \;\; [\,0 \;,\; 1\,]$ , \\
(iv) $a \;\ast\; b \;\, \leq \;\, c \;\ast\; d$
\, whenever \, $a \;\leq\; c$  ,  $b \;\leq\; d$  and  $a \,
, \, b \, , \, c \, , \, d \;\, \in \;\;[\,0 \;,\; 1\,]$.
\end{definition}
A few examples of continuous t-norm are $\,a\,\ast\,b\,=\,ab,\;\,a\,\ast\,b\,=\,min\{a,b\},\;\,a\,\ast\,b\,=\,max\{a+b-1,0\}$.

\begin{definition}
\cite{Schweizer}. A binary operation \, $\diamond \; : \; [\,0 \; ,
\; 1\,] \; \times \; [\,0 \; , \; 1\,] \;\, \rightarrow \;\,
[\,0 \; , \; 1\,]$ \, is continuous \, $t$-conorm if
\,$\diamond$\, satisfies the
following conditions \, $:$ \\
(i)  $\diamond$ \, is commutative and
associative ,\\
(ii)  $\diamond$ \, is continuous , \\
(iii) $a \;\diamond\;0 \;\,=\;\, a
\hspace{1.2cm}
\forall \;\; a \;\; \in\;\; [\,0 \;,\; 1\,]$ , \\
(iv) $a \;\diamond\; b \;\, \leq \;\, c
\;\diamond\; d$ \, whenever \, $a \;\leq\; c$  ,  $b \;\leq\; d$
 and  $a \, , \, b \, , \, c \, , \, d \;\; \in\;\;[\,0
\;,\; 1\,].$
\end{definition}
A few examples of continuous t-conorm are $\,a\,\diamond\,b\,=\,a+b-ab,\;\,a\,\diamond\,b\,=\,max\{a,b\},\;\,a\,\diamond\,b\,=\,min\{a+b,1\}$.\\

We further assume that
$\;(\bf C1)\;$ $a\ast a=a,\;\;\;\;(\bf C2)$ $\;a\diamond a=a.$

\section{Union and intersection of intuitionistic fuzzy soft sets redefined}
\begin{definition}
The intersection of two intuitionistic fuzzy soft sets $(\mathcal{F},A)$ and $(\mathcal{G},B)$ is denoted by $(\mathcal{F},A)\tilde{\cap}\,(\mathcal{G},B)$ and defined by a intuitionistic fuzzy soft set $\mathcal{H}:A\cap B\rightarrow IF^{U}$ such that for each $e\in A\cap B$
\[\mathcal{H}(e)=\{(x,\,\mu_{_{\mathcal{H}(e)}}(x),\,\nu_{_{\mathcal{H}(e)}}(x)):x\in U\}\]
where $\;\;\mu_{_{\mathcal{H}(e)}}(x)=\mu_{_{\mathcal{F}(e)}}(x)\ast \mu_{_{\mathcal{G}(e)}}(x),\;\;\; \nu_{_{\mathcal{H}(e)}}(x)=\nu_{_{\mathcal{F}(e)}}(x)\diamond \nu_{_{\mathcal{G}(e)}}(x).$
\end{definition}

\begin{definition}
The union of two intuitionistic fuzzy soft sets $\mathcal{F}_\alpha$ and $\mathcal{G}_\beta$ is denoted by $(\mathcal{F},A) \tilde{\cup}\,(\mathcal{G},B)$ and defined by a intuitionistic fuzzy soft set $\mathcal{H}:A\cup B\rightarrow IF^{U}$ such that for each $e\in A\cup B$
\[\mathcal{H} (e)=\{(x,\,\mu_{_{\mathcal{F}(e)}}(x),\,\nu_{_{\mathcal{F}(e)}}(x)):x\in U\}\;\;\;if \;e\in A-B\]
\[=\{(x,\,\mu_{_{\mathcal{G}(e)}}(x),\,\nu_{_{\mathcal{G}(e)}}(x)):x\in U\}\;\;\;if \;e\in B-A\hspace{-1 cm}\]
\[=\{(x,\,\mu_{_{\mathcal{H}(e)}}(x),\,\nu_{_{\mathcal{H}(e)}}(x)):x\in U\}\;\;\;if \;e\in A\cap B\hspace{-1 cm}\]
where $\;\;\mu_{_{\mathcal{H}(e)}}(x)=\mu_{_{\mathcal{F}(e)}}(x)\diamond \mu_{_{\mathcal{G}(e)}}(x),\;\;\; \nu_{_{\mathcal{H}(e)}}(x)=\nu_{_{\mathcal{F}(e)}}(x)\ast \nu_{_{\mathcal{G}(e)}}(x).$
\end{definition}

\begin{example}\label{e1}
Suppose that $U$ is the set of medicine under consideration and $U=\{m_1,m_2,m_3,m_4\}$. \\
Let $A=\{$fever, chest pain, cough $\}$ and $B=\{$ fever, cough$\}$. i,e.,  $A=\{f,p,c\}$ and $B=\{f,c\}$.\\
The intuitionistic fuzzy soft set $(\mathcal{F},A)$ be defined as follows:\\
$\mathcal{F}(f)=\left\{(m_1,0.9,0.05),(m_2,0.25,0.6),(m_3,0.65, 0.2), (m_4,0.8,0.1)\right\}$\\
$\mathcal{F}(p)=\left\{(m_1,0.3,0.6),(m_2,0.9,0.1),(m_3,0.4,0.6), (m_4,0.3,0.65)\right\}$\\
$\mathcal{F}(c)=\left\{(m_1,0.6,0.2),(m_2,0.3,0.6),(m_3,0.9,0.05), (m_4,0.85,0.1)\right\}$.\\
The intuitionistic fuzzy soft set $(\mathcal{G},B)$ be defined as follows:\\
$\mathcal{G}(f)=\left\{(m_1,0.85, 0.1),(m_2,0.2,0.7),(m_3,0.5,0.4), (m_4,0.8,0.1)\right\}$\\
$\mathcal{G}(c)=\left\{(m_1,0.65, 0.3),(m_2,0.3,0.65),(m_3,0.9,0.1), (m_4,0.7,0.2)\right\}$.\\
Let us consider the t-norm function $\ast$ and t-conorm function $\diamond$ as $a\ast b=ab$ and $a\diamond b=a+b-ab.$ Then\\
$(\mathcal{F}\tilde{\cap}\mathcal{G})(f)=\left\{(m_1,0.765,0.145),(m_2,0.05,0.88),(m_3,0.325,0.52), (m_4,0.64,0.19)\right\}$\\
$(\mathcal{F}\tilde{\cap}\mathcal{G})(c)=\left\{(m_1,0.39,0.44),(m_2,0.09,0.86),(m_3,0.81,0.145), (m_4,0.595,0.28)\right\}$\\
$(\mathcal{F}\tilde{\cup}\mathcal{G})(f)=\left\{(m_1,0.985,0.005),(m_2,0.4,0.42),(m_3,0.825,0.08), (m_4,0.96,0.01)\right\}$\\
$(\mathcal{F}\tilde{\cup}\mathcal{G})(p)=\mathcal{F}(p)$\\
$(\mathcal{F}\tilde{\cup}\mathcal{G})(c)=\left\{(m_1,0.86, 0.06),(m_2,0.51,0.39),(m_3,0.99,0.005), (m_4,0.955,0.02)\right\}$.
\end{example}

\begin{theorem}
Let $(\mathcal{F},A)$, $(\mathcal{G},B)$ and $(\mathcal{H},C)$ be any three intuitionistic fuzzy soft sets over $(U,E),$ then the following holds:\\
(i) $(\mathcal{F},A)\tilde{\cup}\,(\mathcal{G},B)=(\mathcal{G},B)\tilde{\cup}\,(\mathcal{F},A).\\$
(ii) $(\mathcal{F},A)\tilde{\cap}\,(\mathcal{G},B)=(\mathcal{G},B)\tilde{\cap}\,(\mathcal{F},A).\\$
(iii) $(\mathcal{F},A)\tilde{\cup}\,((\mathcal{G},B) \tilde{\cup}\,(\mathcal{H},C))=((\mathcal{F},A)\tilde{\cup}\,(\mathcal{G},B))\tilde{\cup}\,(\mathcal{H},C) .\\$
(iv) $(\mathcal{F},A)\tilde{\cap}\, ((\mathcal{G},B)\tilde{\cap}\,(\mathcal{H},C))=((\mathcal{F},A)\tilde{\cap}\,(\mathcal{G},B))\tilde{\cap}\,(\mathcal{H},C).$
\end{theorem}
{\bf Proof.} Since the $t$-norm function and $t$-conorm functions are commutative and associative, therefore the theorem follows.

\begin{remark}
Let $(\mathcal{F},A)$, $(\mathcal{G},B)$ and $(\mathcal{H},C)$ be any three intuitionistic fuzzy soft sets over $(U,E)$. If we consider $a\ast b=\min\{a,\,b\}$ and $a\diamond b=\max\{a,\,b\}$ then the following holds:\\
(i) $(\mathcal{F},A)\tilde{\cap}\,((\mathcal{G},B) \tilde{\cup}\,(\mathcal{H},C))=((\mathcal{F},A)\tilde{\cap}\,(\mathcal{G},B))\tilde{\cup}\,((\mathcal{F},A)\tilde{\cap}\,(\mathcal{H},C)). \\$
(ii) $(\mathcal{F},A)\tilde{\cup}\,((\mathcal{G},B) \tilde{\cap}\,(\mathcal{H},C))=((\mathcal{F},A)\tilde{\cup}\,(\mathcal{G},B))\tilde{\cap}\,((\mathcal{F},A)\tilde{\cup}\,(\mathcal{H},C)) .\\$
But in general above relations do not hold.
\end{remark}

\section{Relations on intuitionistic fuzzy soft sets}
\begin{definition}
Let $U$ be an initial universal set and $E$ be the set of parameters. Let $A,B\subseteq E$ and $(\mathcal{F},A),\,(\mathcal{G},B)$ be two intuitionistic fuzzy soft sets over $(U,E)$. Then the cartesian product of  $(\mathcal{F},A)$ and $(\mathcal{G},B)$ is denoted by  $(\mathcal{F},A)\times(\mathcal{G},B)=(\mathcal{H},C)$ where $C=A\times B$ and $\mathcal{H}:C\rightarrow IF^U$ is defined as
\[\mathcal{H}(a,b)=(\mathcal{F},A)\tilde{\cap}(\mathcal{G},B),\;\;\;(a,b)\in C.\]
\end{definition}

\begin{example}
Consider the example \ref{e1}. Let $(\mathcal{F},A)\times(\mathcal{G},B)=(\mathcal{H},C)$, where $C=A\times B$. Consider the t-norm function $\ast$ and t-conorm function $\diamond$ as $a\ast b=ab$ and $a\diamond b=a+b-ab.$ Then the elements of $(\mathcal{H},C)$ will be as follows:\\
$\mathcal{H}(f,f)=\left\{(m_1,0.765,0.145),(m_2,0.05,0.88),(m_3,0.325,0.52), (m_4,0.64,0.19)\right\}$\\
$\mathcal{H}(f,c)=\left\{(m_1,0.585,0.335),(m_2,0.075,0.86),(m_3,0.585,0.28), (m_4,0.56,0.28)\right\}$\\
$\mathcal{H}(p,f)=\left\{(m_1,0.255,0.64),(m_2,0.18,0.73),(m_3,0.2,0.76), (m_4,0.24,0.685)\right\}$\\
$\mathcal{H}(p,c)=\left\{(m_1,0.195,0.72),(m_2,0.27,0.685),(m_3,0.36,0.64), (m_4,0.21,0.72)\right\}$\\
$\mathcal{H}(c,f)=\left\{(m_1,0.51,0.28),(m_2,0.06,0.88),(m_3,0.45,0.43), (m_4,0.68,0.19)\right\}$\\
$\mathcal{H}(c,c)=\left\{(m_1,0.39,0.44),(m_2,0.09,0.86),(m_3,0.81,0.145), (m_4,0.595,0.28)\right\}$.
\end{example}

\begin{definition}
Let $(\mathcal{F},A),\,(\mathcal{G},B)$ be two intuitionistic fuzzy soft sets over $(U,E)$. Then an {\bf intuitionistic fuzzy soft relation} from  $(\mathcal{F},A)$ to $(\mathcal{G},B)$ is an intuitionistic fuzzy soft subset of  $(\mathcal{F},A)\times(\mathcal{G},B)$.\\
In other words, an intuitionistic fuzzy soft relation from  $(\mathcal{F},A)$ to $(\mathcal{G},B)$ is of the form $(\mathcal{R},C)$ where $C\subseteq A\times B$ and $\mathcal{R}(a,b)\subseteq(\mathcal{F},A)\times(\mathcal{G},B),\;\;\forall (a,b)\in C$.
\end{definition}

\begin{example}
Consider the intuitionistic fuzzy soft sets $(\mathcal{F},A)$ and $(\mathcal{G},B)$ over $U$ defined as follows:\\
Let $U$ be the set of medicines under consideration and $U=\{m_1,m_2,m_3,m_4\}$.\\
$A$ and $B$ describes the disease of the patients\\
$A=\{$maleria, dengue$\}$ and $B=\{$fileria, cough$\}$. i.e.,  $A=\{m,d\}$ and $B=\{f,c\}$.\\
Define an intuitionistic fuzzy soft relation $(\mathcal{R},C)$ from $(\mathcal{F},A)$ to $(\mathcal{G},B)$ as follows:\\
$(a,b)\in C\subseteq A\times B$ if and only if $a$ and $b$ are both from mosquito.\\
Then the intuitionistic fuzzy soft relation $(\mathcal{R},C)=\{\mathcal{F}(m)\times\mathcal{G}(f),\mathcal{F}(d)\times\mathcal{G}(f)\}$.
\end{example}

\begin{example}
Let $U$ be the set of candidates attending for an interview of school service and $U=\{c_1,c_2,c_3,c_4,c_5,c_6,c_7\}$.\\
$A$ denotes the academic qualification of candidates\\
$A=\{$B.Sc., M.Com., M.A., M.Sc$\}$  i.e., $A=\{b,c,a,s\}$.\\
Then the intuitionistic fuzzy soft set $(\mathcal{F},A)$ describes the B.Sc. passed students, M.Com. passed students and so on. \\
A relation $(\mathcal{R},C)$ on $(\mathcal{F},A)$ be defined by\\
$(a,b)\in C\subseteq A\times A$ if and only if $a$ and $b$ both are of science.\\
Then $(\mathcal{R},C)=\{\mathcal{F}(b)\times\mathcal{F}(b),\mathcal{F}(b)\times\mathcal{F}(s), \mathcal{F}(s)\times\mathcal{F}(s),\mathcal{F}(s)\times\mathcal{F}(b)\}$.
\end{example}

\begin{definition}
Let $\mathcal{R}$ be an intuitionistic fuzzy soft relation from  $(\mathcal{F},A)$ to $(\mathcal{G},B)$ then $\mathcal{R}^{-1}$ is defined as
\[\mathcal{R}^{-1}(a,b)=\mathcal{R}(b,a),\;\;\;\forall (a,b)\in C\subseteq A\times B\]
\end{definition}

\begin{proposition}
If $\mathcal{R}$ is an intuitionistic fuzzy soft relation from  $(\mathcal{F},A)$ to $(\mathcal{G},B)$ then $\mathcal{R}^{-1}$ is a intuitionistic fuzzy soft relation from  $(\mathcal{G},B)$ to $(\mathcal{F},A)$.
\end{proposition}
{\bf Proof.} $\mathcal{R}^{-1}(a,b)=\mathcal{R}(b,a)=G(b)\tilde{\cap}F(a),\;\;\;\;\forall a,b\in C\subseteq A\times B.\\$
Hence $\mathcal{R}^{-1}$ is an intuitionistic fuzzy soft relation from  $(\mathcal{G},B)$ to $(\mathcal{F},A)$.

\begin{proposition}
If $\mathcal{R}_1$ and $\mathcal{R}_2$ be two intuitionistic fuzzy soft relations from $(\mathcal{F},A)$ to $(\mathcal{G},B)$ then\\
(i) $({\mathcal{R}_1}^{-1})^{-1}=\mathcal{R}_1$ and (ii) $\mathcal{R}_1\subseteq \mathcal{R}_2\Rightarrow {\mathcal{R}_1}^{-1}\subseteq {\mathcal{R}_2}^{-1}.$
\end{proposition}
{\bf Proof.} (i) $({\mathcal{R}_1}^{-1})^{-1}(a,b)={\mathcal{R}_1}^{-1}(b,a)=\mathcal{R}_1(a,b)$. Hence $({\mathcal{R}_1}^{-1})^{-1}=\mathcal{R}_1$.\\
(ii) $\mathcal{R}_1(a,b)\subseteq \mathcal{R}_2(a,b)\Rightarrow{\mathcal{R}_1}^{-1}(b,a)\subseteq {\mathcal{R}_2}^{-1}(b,a)\Rightarrow {\mathcal{R}_1}^{-1}\subseteq {\mathcal{R}_2}^{-1}.$

\begin{definition}
The composition $\circ$ of two intuitionistic fuzzy soft relations $\mathcal{R}_1$ and $\mathcal{R}_2$ is defined by \\
\[(\mathcal{R}_1\circ \mathcal{R}_2)(a,c)=\mathcal{R}_1(a,b)\tilde{\cap}\mathcal{R}_2(b,c)\]
where $\mathcal{R}_1$ is a intuitionistic fuzzy soft relation form $(\mathcal{F},A)$ to $(\mathcal{G},B)$ and $\mathcal{R}_2$ is a intuitionistic fuzzy soft relation from $(\mathcal{G},B)$ to $(\mathcal{H},C)$.
\end{definition}

\begin{proposition}
If $\mathcal{R}_1$ and $\mathcal{R}_2$ be two intuitionistic fuzzy soft relations from $(\mathcal{F},A)$ to $(\mathcal{G},B)$ then $(\mathcal{R}_1\circ \mathcal{R}_2)^{-1}={\mathcal{R}_2}^{-1}\circ{\mathcal{R}_1}^{-1}$.
\end{proposition}
{\bf Proof.} $(\mathcal{R}_1\circ \mathcal{R}_2)^{-1}(a,c)=(\mathcal{R}_1\circ \mathcal{R}_2)(c,a)=\mathcal{R}_1(c,b)\tilde{\cap}\mathcal{R}_2(b,a)=\mathcal{R}_2(b,a)\tilde{\cap}\mathcal{R}_1(c,b) ={\mathcal{R}_2}^{-1}(a,b)\tilde{\cap}{\mathcal{R}_1}^{-1}(b,c)= ({\mathcal{R}_2}^{-1}\circ{\mathcal{R}_1}^{-1})(a,c)$.\\
Hence $(\mathcal{R}_1\circ \mathcal{R}_2)^{-1}={\mathcal{R}_2}^{-1}\circ{\mathcal{R}_1}^{-1}$.

\begin{theorem}
$\mathcal{R}_1$ is an intuitionistic fuzzy soft relation form $(\mathcal{F},A)$ to $(\mathcal{G},B)$ satisfying $(C1)$ and $(C2)$ and $\mathcal{R}_2$ is a intuitionistic fuzzy soft relation from $(\mathcal{G},B)$ to $(\mathcal{H},C)$ satisfying $(C1)$ and $(C2)$ then $\mathcal{R}_1\circ \mathcal{R}_2$ is an intuitionistic fuzzy soft relation from $(\mathcal{F},A)$ to $(\mathcal{H},C)$.
\end{theorem}
{\bf Proof.} By definition\\
$\mathcal{R}_1(a,b)= \mathcal{F}(a)\tilde{\cap}\mathcal{G}(b)=
\{\left(x,\,\mu_{_{\mathcal{F}(a)}}(x)\ast \mu_{_{\mathcal{G}(b)}}(x), \,\nu_{_{\mathcal{F}(a)}}(x)\diamond \nu_{_{\mathcal{G}(b)}}(x)\right):x\in U \},\;\forall (a,b)\subseteq A\times B .\\$
$\mathcal{R}_2(b,c)=\mathcal{G}(b)\tilde{\cap}\mathcal{H}(c)=\{\left(x,\,\mu_{_{\mathcal{G}(b)}}(x)\ast \mu_{_{\mathcal{H}(c)}}(x), \,\nu_{_{\mathcal{G}(b)}}(x)\diamond \nu_{_{\mathcal{H}(c)}}(x)\right):x\in U \},\;\forall (b,c)\subseteq B\times C.$
Therefore,\\
$(\mathcal{R}_1\circ \mathcal{R}_2)(a,c)=\mathcal{R}_1(a,b)\tilde{\cap} \mathcal{R}_2(b,c)\\=\{(x,\;(\mu_{_{\mathcal{F}(a)}}(x)\ast \mu_{_{\mathcal{G}(b)}}(x))\;\ast\;(\mu_{_{\mathcal{G}(b)}}(x)\ast \mu_{_{\mathcal{H}(c)}}(x)),
\;(\nu_{_{\mathcal{F}(a)}}(x)\diamond \nu_{_{\mathcal{G}(b)}}(x))\;\diamond\;(\nu_{_{\mathcal{G}(b)}}(x)\diamond \nu_{_{\mathcal{H}(c)}}(x))\,:\;x\in U \},\;\forall (a,b,c)\subseteq A\times B\times C\,.\;$
\\Now,  $(\mu_{_{\mathcal{F}(a)}}(x)\ast \mu_{_{\mathcal{G}(b)}}(x))\;\ast\;(\mu_{_{\mathcal{G}(b)}}(x)\ast \mu_{_{\mathcal{H}(c)}}(x))\\= \mu_{_{\mathcal{F}(a)}}(x)\ast \mu_{_{\mathcal{G}(b)}}(x)\ast \mu_{_{\mathcal{H}(c)}}(x)\\ \leq \mu_{_{\mathcal{F}(a)}}(x)\ast\, 1\,\ast \mu_{_{\mathcal{H}(c)}}(x)\\= \mu_{_{\mathcal{F}(a)}}(x)\ast \mu_{_{\mathcal{H}(c)}}(x)\,\;$\\ and
$\,(\nu_{_{\mathcal{F}(a)}}(x)\diamond \nu_{_{\mathcal{G}(b)}}(x))\;\diamond\;(\nu_{_{\mathcal{G}(b)}}(x)\diamond \nu_{_{\mathcal{H}(c)}}(x))\\= \nu_{_{\mathcal{F}(a)}}(x)\diamond \nu_{_{\mathcal{G}(b)}}(x)\diamond \nu_{_{\mathcal{H}(c)}}(x)\\ \geq \nu_{_{\mathcal{F}(a)}}(x)\diamond\, 0\,\diamond \nu_{_{\mathcal{H}(c)}}(x)\\= \nu_{_{\mathcal{F}(a)}}(x)\diamond \nu_{_{\mathcal{H}(c)}}(x).$
Hence $\mathcal{R}_1(a,b)\tilde{\cap} \mathcal{R}_2(b,c)\subseteq \mathcal{F}(a)\tilde{\cap}\mathcal{H}(c).\\$
Thus $\;\mathcal{R}_1\circ\mathcal{ R}_2$ is an intuitionistic fuzzy soft relation from $(\mathcal{F},A)$ to $(\mathcal{H},C)$.

\section{Symmetric, transitive and reflexive relations on intuitionistic fuzzy soft sets}
\begin{definition}
An intuitionistic fuzzy soft relation $\mathcal{R}$ on $(\mathcal{F},A)$ is said to be intuitionistic fuzzy soft symmetric ralation if $\;\mathcal{R}(a,b)=\mathcal{R}(b,a),\;\;\;\forall\,a,b\in A.$
\end{definition}

\begin{definition}
An intuitionistic fuzzy soft relation $\mathcal{R}$ on $(\mathcal{F},A)$ is said to be intuitionistic fuzzy soft transitive relation if $\;\mathcal{R}\circ \mathcal{R}\subseteq\mathcal{R}.$
\end{definition}

\begin{definition}
An intuitionistic fuzzy soft relation $\mathcal{R}$ on $(\mathcal{F},A)$ is said to be intuitionistic fuzzy soft reflexive relation if $\;\mathcal{R}(a,b)\subseteq \mathcal{R}(a,a)$ and $\mathcal{R}(b,a)\subseteq \mathcal{R}(a,a),\;\;\;\forall\,a,b\in A.$
\end{definition}

\begin{definition}
An intuitionistic fuzzy soft relation $\mathcal{R}$ on $(\mathcal{F},A)$ is said to be an intuitionistic fuzzy soft equivalence relation if it is symmetric, transitive and reflexive.
\end{definition}

\begin{proposition}
If $\mathcal{R}$ is symmetric if and only if $\mathcal{R}^{-1}$ is so.
\end{proposition}
{\bf Proof.} Let $\mathcal{R}$ is symmetric. Then  $\mathcal{R}^{-1}(a,b)=\mathcal{R}(b,a)=\mathcal{R}(a,b)=\mathcal{R}^{-1}(b,a)$. So, $\mathcal{R}^{-1}$ is symmetric.\\
Conversely, let $\mathcal{R}^{-1}$ is symmetric. Then $\mathcal{R}(a,b)=(\mathcal{R}^{-1})^{-1}(a,b)=\mathcal{R}^{-1}(b,a)=\mathcal{R}^{-1}(a,b)=\mathcal{R}(b,a).$ So, $\mathcal{R}$ is symmetric.

\begin{proposition}
$\mathcal{R}$ is symmetric if and only if $\mathcal{R}=\mathcal{R}^{-1}.$
\end{proposition}
{\bf Proof.} Let $\mathcal{R}$ is symmetric. Then  $\mathcal{R}^{-1}(a,b)=\mathcal{R}(b,a)=\mathcal{R}(a,b)$. So, $\mathcal{R}^{-1}=\mathcal{R}.$\\
Conversely, let $\mathcal{R}^{-1}=\mathcal{R}.$ Then $\mathcal{R}(a,b)=\mathcal{R}^{-1}(a,b)=\mathcal{R}(b,a)$. So, $\mathcal{R}$ is symmetric.

\begin{proposition}
If $\mathcal{R}_1$ and $\mathcal{R}_2$ are symmetric relations on $(\mathcal{F},A)$ then $\mathcal{R}_1\circ \mathcal{R}_2$ is symmetric on $(\mathcal{F},A)$ if and only if $\mathcal{R}_1\circ \mathcal{R}_2=\mathcal{R}_2\circ \mathcal{R}_1.$
\end{proposition}
{\bf Proof.} $\mathcal{R}_1$ and $R_2$ are symmetric implies $\mathcal{R}_1^{-1}=\mathcal{R}_1$ and $\mathcal{R}_2^{-1}=\mathcal{R}_2$. Now $(\mathcal{R}_1\circ\mathcal{R}_2)^{-1}=\mathcal{R}_2^{-1}\circ \mathcal{R}_1^{-1}$. So, $\mathcal{R}_1\circ \mathcal{R}_2$ is symmetric implies  $\mathcal{R}_1\circ \mathcal{R}_2=(\mathcal{R}_1\circ\mathcal{R}_2)^{-1}=\mathcal{R}_2^{-1}\circ \mathcal{R}_1^{-1}=\mathcal{R}_2\circ \mathcal{R}_1.\\$
Conversely, $(\mathcal{R}_1\circ\mathcal{R}_2)^{-1}=\mathcal{R}_2^{-1}\circ \mathcal{R}_1^{-1}=\mathcal{R}_2\circ \mathcal{R}_1=\mathcal{R}_1\circ \mathcal{R}_2.$  So, $\mathcal{R}_1\circ\mathcal{R}_2$ is symmetric.

\begin{corollary}
If $\mathcal{R}$ is symmetric then $\mathcal{R}^n$ is symmetric for all positive integer $n$, where $\mathcal{R}^n$ is $\mathcal{R}\circ\mathcal{R}\circ\cdots\circ\mathcal{R}\;(n times)$.
\end{corollary}

\begin{proposition}
If $\mathcal{R}$ is transitive then $\mathcal{R}^{-1}$ is also transitive.
\end{proposition}
{\bf Proof.} $\mathcal{R}^{-1}(a,b)=\mathcal{R}(b,a)\supseteq (\mathcal{R}\circ\mathcal{R})(b,a)=\mathcal{R}(b,c)\tilde{\cap}\mathcal{R}(c,a) = \mathcal{R}(c,a)\tilde{\cap}\mathcal{R}(b,c)=\mathcal{R}^{-1}(a,c)\tilde{\cap}\mathcal{R}^{-1}(c,b)=(\mathcal{R}^{-1}\circ\mathcal{R}^{-1})(a,b)$.\\
So, $\mathcal{R}^{-1}\circ\mathcal{R}^{-1}\subseteq \mathcal{R}^{-1}$. Hence the result.

\begin{proposition}
If $\mathcal{R}$ is transitive then $\mathcal{R}\circ\mathcal{R}$ is so.
\end{proposition}
{\bf Proof.} $(\mathcal{R}\circ\mathcal{R})(a,b)=\mathcal{R}(a,c)\tilde{\cap}\mathcal{R}(c,b) \supseteq (\mathcal{R}\circ\mathcal{R})(a,c)\tilde{\cap}(\mathcal{R}\circ\mathcal{R})(c,b) =(\mathcal{R}\circ\mathcal{R}\circ\mathcal{R}\circ\mathcal{R})(a,b)$.\\
So, $\mathcal{R}\circ\mathcal{R}\circ\mathcal{R}\circ\mathcal{R}\subseteq \mathcal{R}\circ\mathcal{R}.$ Hence the result.

\begin{proposition}
If $\mathcal{R}$ is reflexive then $\mathcal{R}^{-1}$ is so.
\end{proposition}
{\bf Proof.} $\mathcal{R}^{-1}(a,b)=\mathcal{R}(b,a)\subseteq \mathcal{R}(a,a)=\mathcal{R}^{-1}(a,a).$  and\\
$\mathcal{R}^{-1}(b,a)=\mathcal{R}(a,b)\subseteq \mathcal{R}(a,a)=\mathcal{R}^{-1}(a,a).$  \\
Hence the proof.

\begin{proposition}
If $\mathcal{R}$ is symmetric and transitive then $\mathcal{R}$ is reflexive.
\end{proposition}
{\bf Proof.} $\mathcal{R}(a,a)\supseteq (\mathcal{R}\circ\mathcal{R})(a,a)$, since $\mathcal{R}$ is transitive\\
$=\mathcal{R}(a,b)\tilde{\cap}\mathcal{R}(b,a)=\mathcal{R}(a,b)\tilde{\cap}\mathcal{R}(a,b)$, since $\mathcal{R}$ is symmetric\\
$=\mathcal{R}(a,b).\\$
Similarly we can show that $\mathcal{R}(a,a)\supseteq \mathcal{R}(b,a)$.\\
Hence the proof.

\section{Conclusion}
In this paper the theoretical point of view of intuitionistic fuzzy soft set is discussed. we extend the concept of relation in intuitionistic fuzzy soft set theory context. These are supporting structure for research and development of soft set theory.

\end{document}